\documentclass[psamsfonts, reqno]{amsart}

\usepackage{amscd,amsmath,amssymb,amsfonts,verbatim}
\usepackage{times}
\usepackage[cmtip, all]{xy}
\usepackage{graphicx}
\usepackage[active]{srcltx}
\usepackage{color}
\usepackage{textcomp}

\definecolor{refkey}{gray}{.5}   % graylevel for refs
\definecolor{labelkey}{gray}{.5} % graylevel for labels
\definecolor{Red}{rgb}{1,0,0}

\newcommand{\ol}{\overline}		

\newcommand{\pf}{{\bf Proof : }}

\newcommand{\qedwhite}{\hfill \ensuremath{\Box}} 

\newtheorem{theo}{Theorem}[section]
\newtheorem{prop}[theo]{Proposition}
\newtheorem{lem}[theo]{Lemma}
\newtheorem{cor}[theo]{Corollary}

\theoremstyle{definition}

\newtheorem{rem}[theo]{Remark}

\newtheorem{defi}[theo]{Definition}

\newcommand{\Um}{\mbox{\rm Um}}		\newcommand{\SL}{\mbox{\rm SL}}
\newcommand{\GL}{\mbox{\rm GL}}		
	\newcommand{\E}{\mbox{\rm E}}

\title{Applications of Swan's Bertini theorem to unimodular rows}
\author{Manoj K. Keshari and Sampat Sharma}
\newcommand{\Addresses}{{% additional braces for segregating \footnotesize
  \bigskip
  \footnotesize

\textsc{Manoj K. Keshari, Department of Mathematics, IIT Bombay
            Mumbai 400076, INDIA}\par\nopagebreak
  \textit{E-mail:} Manoj K. Keshari \texttt{<keshari@math.iitb.ac.in>}

\medskip
  
 \textsc{Sampat Sharma, School of Basic Sciences, IIT Mandi, Mandi 175005 (H.P), India}\par\nopagebreak
  \textit{E-mail:} Sampat Sharma \texttt{<sampat@iitmandi.ac.in; sampat.iiserm@gmail.com>}

  \medskip

  }}

\begin{document}
\maketitle
\subjclass 2020 Mathematics Subject Classification:{13C10, 19D45, 19G12}

 \keywords {Keywords:}~ {Unimodular row, nice group structure, Mennicke symbol}

 \begin{abstract}
Let $R$ be an affine algebra of dimension $d\geq 4$ over a perfect field $k$ of char $\neq 2$ and $I$ be an ideal of $R$. Then 
\begin{enumerate}
\item $MS_{d+1}(R)$ is uniquely divisible prime to char $k$ if $R$ is reduced and $k$ is infinite with $c.d.(k)\leq 1$.
\item $\Um_{d+1}(R,I)/\E_{d+1}(R,I)$ has nice group structure if $c.d._2(k)\leq 2$.
\item $\Um_{d}(R,I)/\E_{d}(R,I)$ has nice group structure if $k$ is algebraically closed of char $k\neq 2,3$ and either (i) $k = \overline{\mathbb{F}}_{p}$ or (ii) $R$ is normal.

\end{enumerate}
 \end{abstract}
 
 \vskip0.50in

\section{Introduction} {\it All rings are assumed to be commutative noetherian with unity $1$.}

The study of unimodular rows belongs to a classical branch of algebraic $K$-theory which uses
cohomological thinking to address some geometrically-inspired questions in commutative algebra. Let $R$ be a ring and $v = (v_1, \ldots,  v_n)\in \Um_n{(R)}$ be a unimodular row. Then $P = R^n/vR$ is a stably free module. It is
an interesting question whether $P$ is actually free i.e. whether the associated vector bundle
on the affine scheme $\mbox{spec}(R)$ is trivial.  Methods of $\mathbb{A}^1$-homotopy theory, pioneered by Asok, Fasel and
others, studied unimodular rows using homotopy theory. If $R$ is smooth, then unimodular
rows have an interpretation as $\mathbb{A}^1$-homotopy classes of maps from $\mbox{spec} (R)$ to a motivic sphere
i.e. an algebro-geometric version of cohomotopy sets. In general, these may not have group
structures, just as in topology, but there are conditions under which group structures are available and thus
more amenable to cohomological methods.

The present paper removes the limitation that motivic homotopy theory imposes, namely
that $\mbox{spec} (R)$ must be smooth. Indeed, it is a general principle  in $K$-theory that
$\mathbb{A}^{1}$-homotopy invariance fails away from smooth or regular settings. There are classical tools to
address this, chiefly Swan's Bertini theorem, which has led Murthy and others to
study vector bundle questions without $\mathbb{A}^1$-invariance. In fact these methods are older
than $\mathbb{A}^{1}$-homotopy theory.

  In this article we extend some of the results on unimodular rows known in the smooth case to non-smooth case essentially using Swan's Bertini theorem. If the result in the smooth case holds for dimension $d\geq \lambda$,
then in the non-smooth case it holds for $d\geq \lambda + 1$.

For a ring $R$ of dimension $d$, the Mennicke symbol $\mbox{MS}_{d+1}(R)$ is a group cooked up
from $\Um_{d+1}(R)$ up to some obvious relations. In the smooth setting it was proved by
Fasel \cite{faselbk} that  $\mbox{MS}_{d+1}(R)$ is isomorphic to the cohomology group $H^d(R; K^{M}_{d+1})$ which houses a ``secondary
obstruction" group for splitting of vector bundles of rank $d$. One would like to think that the Mennicke
symbol group also plays a similar role in the non-smooth setting. We will describe our first result which is on Mennicke symbol group.

Let $R$ be a smooth affine algebra of dimension $d\geq 3$ over a perfect field $k$ of char $\neq 2$ with $c.d.(k)\leq 1$. Then Fasel {\cite[Theorem 2.2]{faselbk}} proved that $MS_{d+1}(R)$ is uniquely divisible prime to the char $k$. Fasel uses methods of $\mathbb A^1$-homotopy theory which works in the smooth case only. Using Swan's Bertini, we generalize it for $d\geq 4$ by removing the smoothness assumption (\ref{20}).

\begin{theo}
Let $R$ be a reduced affine algebra of dimension $d\geq 4$ over an infinite perfect field $k$ of char $\neq 2$ with $c.d.(k)\leq 1$. Then $MS_{d+1}(R)$ is uniquely divisible prime to the char $k$.
\end{theo}

Now we will describe our results on nice group structure on elementary orbit space of unimodular vectors.
Let $R$ be a ring of dimension $d$, $I$ be an ideal of $R$ and $n$ be an integer with $2n\geq d+4$. Then 
van der Kallen \cite{vdk2} proved that the orbit space ${\Um}_n(R,I)/{\E}_n(R,I)$ has  an abelian group structure.
Given $v,w\in {\Um}_n(R,I)$, by  {\cite[Lemma 3.2]{vdk3}},  the classes $[v]$ and $[w]$ can be expressed as $[v]=[(x_1,v_2,\ldots,v_n)]$ and $[w]=[(v_1,v_2,\ldots,v_n)]$. Their product is defined as {\cite[Lemma 3.5]{vdk2}}
$$[(x_1,v_2,\ldots,v_n)]\ast[(v_1,v_2,\ldots,v_n)]=[(v_1(x_1+w_1)-1,(x_1+w_1)v_2,v_3,\ldots,v_n)] $$
  where $w_1\in R$ is such that $v_1w_1=1$ modulo the ideal $\langle v_2,\ldots,v_n\rangle$. In particular, we have a group structure on ${\Um}_d(R,I)/{\E}_d(R,I)$ for $d\geq 4$. Note that for $I=R$, we have ${\Um}_n(R,R)={\Um}_n(R)$ and $\E_n(R,R)=\E_n(R)$.

By van der Kallen-Rao {\cite[Theorem 6.1]{kallenrao}} and Vaserstein {\cite[Theorem 5.2 (a)]{7}}, if $R$ is a smooth affine algebra of dimension $3$ over a perfect $C_1$-field $k$ of char $\neq 2,3$, then
${\Um}_3(R)/{\E}_3(R)$ also has same group structure.
We say that the group operations on ${\Um}_n(R,I)/{\E}_n(R,I)$  is {\it{nice}} if it is given by coordinatewise multiplication i.e.
  $$ [(x_1,v_2,\ldots,v_n)]\ast[(v_1,v_2,\ldots,v_n)]=[(x_1v_1,v_2,\ldots,v_n)].$$
The group structure on ${\Um}_{d+1}(R)/{\E}_{d+1}(R)$ is not nice in general, see {\cite[Example 2.2(c)]{weibelci}} for $d=3$ and \cite[Examples 4.14, 4.17]{vdk1} for $d>3$.

Let $R$ be an affine algebra of dimension $d\geq 3$ over a perfect field $k$ of char $\neq 2$. Then the group structure on 
${\Um}_{d+1}(R)/{\E}_{d+1}(R)$ is nice in the following cases.
\begin{enumerate} 
\item
When $R$ is smooth and c.d.$_{2}(k)\leq 2$ due to Fasel {\cite[Theorem 2.1]{faselbk}}. 
\item When c.d.$(k)\leq 1$ due to Garge-Rao  {\cite[Theorem 3.9]{gr}}. 
\end{enumerate}
Using Swan's Bertini, we generalize Fasel's result by removing the smoothness assumption for $d\geq 4$ (\ref{1}). 

\begin{theo}
Let $R$ be an affine algebra of dimension $d\geq 4$ over a perfect field $k$ of char $\neq 2$ with c.d.$_{2}(k)\leq 2$. Let $J$ be an ideal of $R$. Then the group structures on 
${\Um}_{d+1}(R)/{\E}_{d+1}(R)$ and ${\Um}_{d+1}(R, J)/{\E}_{d+1}(R, J)$ are nice.
\end{theo}

Let $R$ be a smooth affine algebra of dimension $d\geq 3$ over an algebraically closed field $k$ of char $\neq 2,3$ and $I$ be an ideal of $R$. Then Gupta-Garge-Rao {\cite[Theorem 6.8]{ggr}} proved that the group structure on ${\Um}_{d}(R)/{\E}_{d}(R)$ is nice. Further when $d\geq 4$,  the group structure on ${\Um}_{d}(R,I)/{\E}_{d}(R,I)$ is also nice {\cite[Theorem 8.2]{ggr}}. Using Swan's Bertini, we generalize these results for $d\geq 4$ by relaxing the smoothness assumption (\ref{2}, \ref{21}, \ref{22}). This result was proved in {\cite[Theorem 3.4]{kesharisharma}} in absolute case when $d\geq 5$ and $k=\overline{\mathbb{F}}_{p}$.

\begin{theo}
\label{mainthm} Let $R$ be an affine algebra of dimension $d\geq 4$ over an algebraically closed field $k$ of char $\neq 2,3$ and $I$ be an ideal of $R$. Assume that either (i) $k = \overline{\mathbb{F}}_{p}$ or (ii) $R$ is normal. Then 
the group structures on ${\Um_{d}(R)}/{\E_{d}(R)}$ and ${\Um_{d}(R,I)}/{\E_{d}(R,I)}$ are nice.
\end{theo}
 
As a consequence of this we prove the following result.
\begin{cor}
            Let $R$ be an affine algebra of dimension $d\geq 4$ over an algebraically closed field $k$ of char $\neq 2, 3.$ Let $I$ be an ideal of $R.$ Assume that either (i) $k = \overline{\mathbb{F}}_{p}$ or (ii) $R$ is normal. If $\sigma \in {\SL}_{d}(R, I)\cap {\E}_{d+1}(R, I)$, then 
          $e_{1}\sigma$ can be completed to a relative elementary matrix $\varepsilon \in \E_{d}(R, I).$
           \end{cor}

Now we will describe our last set of results on first row map and torsion elements.
Let $R$ be a ring of dimension $d\geq 2$. In \cite[Theorem 3.16(iv)]{vdk1},
van der Kallen  proved that the 
first row map 
$row_1:\SL_{d+1}(R)\rightarrow {\Um_{d+1}(R)}/{\E_{d+1}(R)}$
is a group homomorphism. Further, he gave an example \cite[Examples 4.16, 4.13]{vdk1} showing that 
$row_1:\GL_{d+1}(R)\rightarrow {\Um_{d+1}(R)}/{\E_{d+1}(R)}$
is not a group homomorphism. 

We prove the following results
(\ref{elem3}, \ref{3}, \ref{31}). When $R$ is a local ring, (1) and (2) are due to second author \cite[Theorem 1.1]{grouphom} and \cite[Corollary 3.8]{absence} and $(3)$ is proved in {\cite[Corollary 4.5]{kesharisharma}} in case $I=R$.

\begin{theo}
Let $R$ be a ring of dimension $d\geq 2$ and $I$ be an ideal of $R$. Assume that either (i) $R$ is an affine $C$-algebra with either $C=\mathbb Z$ or $C$ is a subfield of $\overline{\mathbb{F}}_{p}$ or (ii) height of the Jacobson radical $ht(J(R))\geq 1$. Then following holds.
\begin{enumerate}
\item The map $row_1:\GL_{d+1}(R[X]) \rightarrow \Um_{d+1}(R[X])/{\E_{d+1}(R[X])}$ is a group homomorphism.
\item Assume $1/m,1/d!\in R$. Then $\Um_{d+1}(R[X])/{\E_{d+1}(R[X])}$ has no $m$-torsion.
\item Assume $1/2\in R$ and $d\geq 4$. Then the group structure on $\Um_{d+1}(R[X],I[X])/{\E}_{d+1}(R[X],I[X])$ is nice.
\end{enumerate}
\end{theo}

\section{Preliminaries}
\par 
\begin{defi}
\label{definitions} Let $R$ be a ring and $I$ be an ideal of $R$. 
\begin{enumerate}
\item The ideal generated by $a_1,\ldots,a_r\in R$ is $\langle a_1,\ldots,a_r\rangle$.

\item A vector 
 $v = (a_{1},\ldots, a_{r})\in R^r$ is {\it unimodular} if there exist $w=(b_{1},\ldots, b_{r})\in R^r$ with 
$\langle v ,w\rangle = \Sigma_{i = 1}^{r} a_{i}b_{i} = 1$. The set 
of unimodular rows in $R^r$ is $\Um_r(R)$.   
The set $\Um_{r}(R,I)$ consist of unimodular vectors $v\in\Um(R)$ which are congruent to $e_r$ modulo $I$, where $e_r=(0,\ldots,0,1)$. 

\item The group $\E_r(R)$ is generated by elementary matrices $e_{ij}(\lambda) = \mbox{Id}_{r} + \lambda {\E}_{ij} \in \SL_r(R)$, where $\lambda \in R$ and for $i\neq j$,  ${\E}_{ij} \in \mbox{M}_{r}(R)$ whose $(i,j)^{th}$ entry is $1$ and all other entries are zero. 

\item The group
${\E}_{r}(R)$ acts on $\Um_r(R)$ by right multiplication i.e.  $v\cdot\varepsilon := v\varepsilon$ and
 ${\Um_{r}(R)}/{\E_{r}(R)} =\{[v]:v\in \Um_r(R)\}$ denote the elementary orbit space.

\item The group $\mbox{E}_{n}(I)$ is the subgroup of $\mbox{E}_{n}(R)$ generated
 by elements $e_{ij}(x)$ with $x\in I.$ The relative elementary group $\mbox{E}_{n}(R,I)$ is the normal closure of $\mbox{E}_{n}(I)$ in $\mbox{E}_{n}(R).$
  
\item The stable range condition $sr_{n}(R)$ holds for $R$ if 
 given $(a_{1},a_{2},\ldots,a_{n+1})\in \Um_{n+1}(R),$ there exist $c_{i}\in R$ such that 
 $(a_{1}+ c_{1}a_{n+1}, a_{2}+c_{2}a_{n+1},\ldots,a_{n}+c_{n}a_{n+1})\in \Um_{n}(R).$
 
\item The stable range $sr(R)$ of
$R$ is the least integer $n$ such that $sr_{n}(R)$ holds. The stable dimension of $R$ is
$\mbox{sdim}(R) = sr(R) - 1.$
 
\item\label{excision} The excision ring $R\oplus I$ has coordinatewise addition. The multiplication
 is given as 
 $(r,i).(s, j) = (rs, rj+si+ij)$, where $r,s \in R$ and $i,j\in I.$
 The multiplicative identity is $(1,0).$ 
\item If $R$ is reduced, then $R\oplus I$ is reduced. To see this, $(r,i)^n=(r^n,j)=0$ implies $r=0$. Now $(0,i)^n=(0,i^n)=0$ implies $i=0$.

\item A Mennicke symbol of length $n\geq 3$ is a pair $(\phi,G)$, where $G$ is a group and $\phi :\Um_n(R)\to G$ is a map satisfying
\begin{enumerate}
\item $\phi(e_1)=1$ and $\phi (v)=\phi(v\sigma)$ for $\sigma\in \E_n(R)$.
\item $\phi(a,a_2,\ldots,a_n)\cdot \phi(b,a_2,\ldots,a_n)=\phi(ab,a_2,\ldots,a_n)$.
\end{enumerate}
A universal Mennicke symbol exist and is denoted as $(ms,MS_n(R))$.

\end{enumerate}
\end{defi}

We will state a part of Swan's Bertini needed in our proof \cite[Theorems 1.3, 1.4]{swanbertini} (see also \cite[Theorem 2.3]{murthy}).

\begin{theo}\label{sb}
Let $R$ be a geometrically reduced affine ring of dimension $d$ over an infinite field $k$ and $(a,a_1,\ldots,a_r)\in \Um_{r+1}(R)$. Then there exist $b_1,\ldots,b_r\in R$ such that if $a'=a+a_1b_1+\ldots+a_rb_r$, then the hypersurface $R/\langle a'\rangle$ is a reduced affine algebra  of dimension $d-1$ which is smooth at all smooth points of $R$. 
\end{theo}

The following result is due to Garge-Gupta-Rao in \cite[Lemma 3.6]{ggr}.

\begin{lem}\label{relnicecri} Let $R\oplus I$ be the excision ring of $R$ with respect to an ideal $I$ in $R$ and $n\geq 3.$ Suppose both $\Um_{n}(R, I)/\E_{n}(R,I)$ and $\Um_{n}(R\oplus I)/\E_{n}(R\oplus I)$ have group structure given by van der Kallen's rule. Then the group structure on $\Um_{n}(R, I)/\E_{n}(R,I)$ is nice whenever it is nice for $\Um_{n}(R\oplus I)/\E_{n}(R\oplus I).$
\end{lem}

The following local-global principle for unimodular rows is due to Rao {\cite[Theorem 2.3]{topcoeffie}}. 
\begin{theo}
\label{raolg} Let $R$ be a ring, $v\in \Um_{n}(R[X])$ and $n\geq 3.$ Suppose that  $v\equiv v(0)~\mbox{mod}~\E_{n}(R_{\mathfrak{m}}[X])$ for all maximal ideal $\mathfrak{m}$ of $R$. Then $v\equiv v(0)~\mbox{mod}~\E_{n}(R[X]).$
\end{theo}

Next we note a pre-stabilization result of van der Kallen {\cite[Theorem 2.2]{vdk2}}. 

\begin{theo}
\label{prestab} Let $R$ be a ring of stable dimension $\leq 2n-3$ with $n\geq 3$ and $I$ be an ideal of $R.$ Let $i,j\geq 0$ and $g\in \GL_{n+i}(R)\cap \E_{n+i+j+1}(R, I).$ Then there exist matrices $u,v,w, M$ with entries in $I$ and 
$q$ with entries in $R$ such that 
$$ \begin{bmatrix}
                 I_{i+1} + uq & v \\
                 wq & I_{n-1}+M\\
                \end{bmatrix}\in \E_{n+i}(R,I) ~\textit{and}~ \begin{bmatrix}
                 I_{j+1} + qu & qv \\
                 w & I_{n-1}+M\\
                \end{bmatrix}\in g\E_{n+j}(R,I).$$

\end{theo}

Next we prove a key lemma.

\begin{lem}
\label{keylemma} Let $k$ be a field and $d, r(d)\geq 3$ be integers. Let $\Omega_{d}$ be the set of affine $k$-algebras of dimension $d.$ Assume that for all smooth $A\in \Omega_d,$ the group structure on $\Um_{r(d)}(A)/\E_{r(d)}(A)$ is nice. Let 
$R\in \Omega_{d+1}$ and $J$ be the ideal defining singular locus of $R$. Given $v, w\in \Um_{r(d)+1}(R,J)$ with $v = (a, a_1,\ldots, a_{r(d)}), w = (b, a_1,\ldots, a_{r(d)}),$ we have $[v]\ast[w] = [ab, a_1, \ldots, a_{r(d)}].$
\end{lem}

${\pf}$ By Swan's Bertini (\ref{sb}), we can add a linear combination of $ab, a_{1}, \ldots, a_{r(d)-1}$ to $a_{r(d)}$ to get 
$a_{r(d)}' = a_{r(d)} + \lambda ab + \sum_{i = 1}^{r(d)-1}\lambda_{i}a_{i}$ for some $\lambda, \lambda_{i} \in R$ such that ${R}/{\langle a_{r(d)}'\rangle}$ is smooth outside the singular set of $R.$ Since $a_{r(d)}' \equiv 1~(\mbox{mod}~J)$, the hypersurface $\overline R ={R}/{\langle a_{r(d)}'\rangle}$ 
is smooth of dimension $d.$
We have 
$$[v] = [(a, a_{1}, \ldots,a_{r(d)-1}, a_{r(d)}')], ~~~~[w] = [(b, a_{1}, \ldots, a_{r(d)-1}, a_{r(d)}')]$$
 Choose $p\in R$ such that $ap \equiv 1 ~\mbox{mod}~\langle a_1,\ldots, a_{r(d)}'\rangle.$ Then 
$$[v] \ast [w] = [(a(b+p)-1, a_{1}(b+p), a_{2}, \ldots, a_{r(d)}')]$$ 
Since $\overline{R}\in \Omega_d $ is smooth,
\begin{eqnarray*}
[(\overline{a}(\overline{b}+\overline{p})-1, \overline{a_{1}}(\overline{b}+\overline{p}), \overline{a_{2}}, \ldots, \overline{a_{r(d)-1}})]  &=& [(\overline{a}, \overline{a_{1}}, \ldots, \overline{a_{r(d)-1}})]\ast [(\overline{b}, \overline{a_{1}}, \ldots, \overline{a_{r(d)-1}})] \\
&=& [(\overline{ab}, \overline{a_{1}}, \ldots, \overline{a_{r(d)-1}})]
\end{eqnarray*}
Upon lifting the elementary matrix, we have 
\begin{eqnarray*}
[(a, a_{1}, \ldots, a_{r(d)}')]\ast  [(b, a_{1}, \ldots, a_{r(d)}')] &=& [(a(b+p)-1, a_{1}(b+p), a_{2}, \ldots, a_{r(d)}')] \\
&=&  [(ab, a_{1}, \ldots, a_{r(d)}')]
\end{eqnarray*}
 Thus we have 
$$ [v]*[w]= [(a, a_{1}, \ldots, a_{r(d)})]\ast  [(b, a_{1}, \ldots, a_{r(d)})] =  [(ab, a_{1}, \ldots, a_{r(d)})].$$
$\hfill\square$

%%%%%%%%%%%%%%%%%%%%%%%%%%%%%%%%%%%%%%%%%%%

\section{Divisibility of $MS_{d+1}(R)$}

\begin{theo}\label{20}
Let $R$ be a reduced affine algebra of dimension $d\geq 4$ over an infinite perfect field $k$ of char $\neq 2$ with $c.d.(k)\leq 1$. Then $MS_{d+1}(R)$ is uniquely divisible prime to the characteristic of $k$.
\end{theo}

${\pf}$ Let $n\geq 1$ be prime to the characteristic of $k$. If $I$ is the ideal defining the singular locus of $R$, then $ht(I)\geq 1$.  Let $[v]=[(a_1,\ldots,a_{d+1})] \in MS_{d+1}(R)$. Since $dim (R/I)\leq d-1$, we may assume that $v=e_1$ mod $I$.
By Swan's Bertini (\ref{sb}), we can add multiples of $a_2,\ldots,a_{d+1}$ to $a_1$ and assume that the hypersurface $\ol R=R/\langle a_1\rangle$ is of dimension $d-1\geq 3$ and is smooth outside the singular set of $R$. Since $a_1=1$ mod $I$, we get $\ol R$ is smooth affine $k$-algebra of dimension $d-1$. Applying Fasel's result in the smooth case 
{\cite[Theorem 2.2]{faselbk}}, we get that $[\ol v]=[(\ol a_2,\ldots,\ol a_{d+1})] \in MS_d(\ol R)$ is uniquely divisible by $n$. So there exist a unique $[\ol w]=[(\ol b_2,\ldots,\ol b_{d+1})] \in MS_d(\ol R)$ such that $[\ol v]=[\ol w]^n = [(\ol b_2^n,\ol b_3,\ldots,\ol b_{d+1})]\in MS_d(\ol R)$. Therefore, $[v]=[(a_1,b_2^n,\ldots,b_{d+1})]=[(a_1,b_2,\ldots,b_{d+1})]^n \in MS_{d+1}(R)$. Thus $[v]$ is divisible by $n$ in $MS_{d+1}(R)$.

To prove uniqueness, assume there exist $[w_1]=[(c,c_2,\ldots,c_{d+1})]$ and $[w_2]=[(d,c_2,\ldots,c_{d+1})]\in MS_{d+1}(R)$ such that $[v]=[w_1]^n=[w_2]^n$. In view of \cite[Proposition 3.8]{gr}, we may assume that $w_1,w_2=e_1$ mod $I$. By Swan's Bertini (\ref{sb}), adding multiples of $cd,c_2,\ldots,c_d$ to $c_{d+1}$, we may assume that $\ol R = R/\langle c_{d+1}\rangle$ is smooth affine of dimension $d-1$. Since  $[\ol v]=[\ol c,\ol c_2,\ldots,\ol c_d]^n=[\ol d,\ol c_2,\ldots,\ol c_d]^n \in MS_d(\ol R)$, by uniqueness part in Fasel {\cite[Theorem 2.2]{faselbk}}, we get $[\ol c,\ol c_2,\ldots,\ol c_d]=[\ol d,\ol c_2,\ldots,\ol c_d] \in MS_d(\ol R)$. Therefore, lifting elementary automorphism, we have $$[ c + a_1c_{d+1}, c_2 + a_2c_{d+1},\ldots, c_d+a_dc_{d+1}, c_{d+1}]=[d+b_1c_{d+1}, c_2 +b_{2}c_{d+1},\ldots, c_d+b_{d}c_{d+1}, c_{d+1}] $$ for some $a_i, b_i \in R.$ Therefore upon making suitable elementary transformations, $[w_1]=[w_2]\in MS_{d+1}(R)$.
$\hfill \square$

\section{A nice group structure on $\Um_{d+1}(R,I)/\E_{d+1}(R,I)$}

\begin{theo}\label{1}
Let $R$ be an affine algebra of dimension $d\geq 4$ over a perfect field $k$ of char $\neq 2$ with c.d.$_{2}(k)\leq 2$.  Let $J$ be an ideal of $R$. Then the group structures on ${\Um}_{d+1}(R)/{\E}_{d+1}(R)$ 
and ${\Um}_{d+1}(R, J)/{\E}_{d+1}(R, J)$ are nice.
\end{theo}

${\pf}$  We first prove the absolute case ${\Um}_{d+1}(R)/{\E}_{d+1}(R)$. When $k$ is a finite field, the result follows from {\cite[Theorem 3.9]{gr}}. Hence we assume that $k$ is infinite. Further by \cite[Lemma 3.5]{gr}, we may assume $R$ is reduced. Since $k$ is perfect, we get $R$ is geometrically reduced.

{\it{Case I:}} If $I$ denotes the ideal defining the singular locus of $R,$ then ht$(I)\geq 1$. Let $v, w \in \Um_{d+1}(R)$. In view of \cite[Proposition 3.8]{gr}, we may assume that 
$v = (a, a_{1}, \ldots, a_{d})$ and $w = (b, a_{1}, \ldots, a_{d})$ with $v,w\equiv e_{d+1}$ $\text{mod}~I$. In view of Fasel {\cite[Theorem 2.1]{faselbk}} which proves the result for $d\geq 3$ whenever $R$ is smooth and applying (\ref{keylemma}) with $r(d) =d$, we get 
$$ [v]*[w]= [(a, a_{1}, \ldots, a_{d})]\ast  [(b, a_{1}, \ldots, a_{d})] =  [(ab, a_{1}, \ldots, a_{d})].$$

{\it{Case II :}} By {\cite[Proposition 3.1]{K}}, $R\oplus J$ is an affine $k$-algebra of dimension $d$. Since $R$ is a reduced ring, the excision ring $R\oplus J$ is also reduced (\ref{excision}). By case I,  ${\Um}_{d+1}(R\oplus I)/{\E}_{d+1}(R\oplus I)$ has nice group structure. Therefore by  (\ref{relnicecri}), the group structure on ${\Um}_{d+1}(R, J)/{\E}_{d+1}(R, J)$ is nice.
$\hfill \square$

\begin{cor}
           \label{vandercor1}
Let $R$ be an affine algebra of dimension $d\geq 4$ over a perfect field $k$ of char $\neq 2$ with c.d.$_{2}(k)\leq 2$. Let $I$ be an ideal of $R$.
   If $\sigma \in {\SL}_{d+1}(R, I)\cap {\E}_{d+2}(R, I)$, then 
          $e_{1}\sigma$ can be completed to a relative elementary matrix $\varepsilon \in \E_{d+1}(R, I).$
           \end{cor}
    ${\pf}$  In view of  theorem \ref{prestab}, we have
$$ [\sigma] = \begin{bmatrix}
                 1 + ux & uy \\
                 z^{t} & M\\
                \end{bmatrix},
                ~~ [\varepsilon] = \begin{bmatrix}
                 1 + ux & y \\
                 uz^{t} & M\\
                \end{bmatrix},$$ for some $u,x\in R, y,z\in M_{1,d}(R), M\in M_{d,d}(R), 
                \varepsilon \in E_{d+1}(R)$. Since the group structure on ${\Um_{d+1}(R)}/{\E_{d+1}(R)}$ is nice, we have 
$$[e_{1}\sigma] = [1+ux, uy] = [1+ux, u]\ast [1+ux,y] = [1+ux, y] =[e_{1}].$$
The last equality $ [1+ux, y] =[e_{1}]$ follows from the existance of $\varepsilon$.
$\hfill \square$

%%%%%%%%%%%%%%%%%%%
\section{A nice group structure on $\Um_{d}(R,I)/\E_{d}(R,I)$}

Given $v, w\in \Um_{r+1}(R)$, the pair $(v,w)$ is called an {\it{admissible pair}} if $v = (a, a_1, \ldots, a_r)$ and $w = (b, a_1, \ldots, a_r)$. In this case, $c(v,w) = (ab, a_1, \ldots, a_r)\in \Um_{r+1}(R)$  and the triple $([v],[w],[c(v,w)])$ is called the {\it nub} to admissible pair $(v,w).$
A sequence $(v_1, w_1)\rightarrow \cdots \rightarrow (v_n, w_n)$ of admissible pairs is called an {\it admissible sequence} if each pair 
$(v_i,w_i)$ has the same nub. We say that $(v_n, w_n)$ is in the admissible orbit of $(v_1, w_1)$.

The proof of next result is similar to {\cite[Proposition 3.8]{gr}}. 

\begin{prop}
\label{mainprop}Let $R$ be a reduced affine algebra of dimension $d\geq 4$ over an algebraically closed field $k.$ Assume that either (i) $k = \overline{\mathbb{F}}_{p}$ or (ii) $R$ is a normal. Let $I$ be the ideal defining singular locus of $R$ and  $(v, w)$ be an admissible pair of size $d$ over $R.$ Then there exists another admissible pair 
$(v', w')$ of size $d$ over $R$ such that both $v', w' \equiv e_{d}~\mbox{mod}~I$ and $(v, w)\rightarrow (v', w')$ is an admissible sequence.
\end{prop}

${\pf}$  
 {\it{Case 1:}} $k = \overline{\mathbb{F}}_{p}.$ Let $v = (a, a_{2}, \ldots, a_{d}), w = (b, a_{2}, \ldots, a_{d})\in \Um_{d}(R).$ Since $R$ is a reduced affine algebra, we have $\mbox{ht}(I)\geq 1.$ Thus $\overline{R} = R/I$ has  dimension $\leq d-1$ and $d-1 \geq 3.$ By \cite[Corollary 17.3]{7}, $\mbox{sr}(\overline{R})\leq d-1.$  
Thus there exist $\overline{\lambda_{2}}, \ldots, \overline{\lambda_{d}}\in \overline{R}$ such that if $ \overline{a_{i}'} = \overline{a_{i}} + \overline{\lambda_{i}ab}$ for $i=2,\ldots,d$, then
$(\overline{a_{2}'}, \ldots, \overline{a_{d}'})\in \Um_{d-1}(\overline{R})$. Thus we can modify $\overline{v}, \overline{w}$ to $\overline{v_{1}} = (\overline{a}, \overline{a_{2}'}, \ldots, \overline{a_{d}'}), \overline{w_{1}} = (\overline{b}, \overline{a_{2}'}, \ldots, \overline{a_{d}'})$.
Since $(\overline{a_{2}'}, \ldots, \overline{a_{d}'})\in \Um_{d-1}(\overline{R}),$ after making suitable admissible transformations, we may assume that $\overline{a} = \overline{1}$ and $\overline{b} = \overline{1}.$ Next we transform $\overline{v_1}, \overline{w_1}$ to  $\overline{v_{2}} = (\overline{a}, \overline{a_{2}''}, \ldots, \overline{a_{d}''}), \overline{w_{2}} = (\overline{b}, \overline{a_{2}''}, \ldots, \overline{a_{d}''})$, where $\overline{a_{d}''} = \overline{a_{d}'} + \overline{ab}(1-a_{d}')$ and 
$\overline{a_{j}''} = \overline{a_{j}'} + \overline{ab}(-a_{j}')$ for $j=2,\ldots, d-1.$ Thus we have 
$\overline{a_{d}''} = \overline{1}$ and $\overline{a_{j}''} = \overline{0}.$ Further by modifying $\overline{v_2}, \overline{w_2}$ to  $\overline{v_{3}} = (\overline{a} - \overline{a_{d}''}, \overline{a_{2}''}, \ldots, \overline{a_{d}''}), \overline{w_{3}} = (\overline{b} - \overline{a_{d}''}, \overline{a_{2}''}, \ldots, \overline{a_{d}''})$, we may assume that $\overline{a}, \overline{b} = \overline{0}.$  Thus upon lifting the elementary transformations, we can modify $v, w$ to $v_{3} =  ({a}, {a_{2}''}, \ldots, a_{d-1}'', {a_{d}''}), w_{3} = ({b}, {a_{2}''}, \ldots, a_{d-1}'', {a_{d}''})$ respectively such that $v_{3}, w_{3} \equiv e_{d}~\mbox{mod}~I$ and $(v,w)\rightarrow (v_{1}, w_{1})\rightarrow (v_{2}, w_{2}) \rightarrow (v_{3}, w_{3})$ is an admissible sequence.

{\it {Case 2:}} Assume $R$ is normal, then $\mbox{ht}(I)\geq 2.$  Thus $\mbox{dim}(\overline{R})\leq d-2.$ Since $(\overline{ab}, \overline{a_{2}}, \ldots, \overline{a_{d}})\in \Um_{d}(\overline{R}),$ by 
\cite[Corollary 9.4]{7}, there exist $\overline{c_{2}}, \ldots, \overline{c_{d}}\in \overline{R}$ such that $\mbox{ht}(\sum_{i=2}^{d}\overline{R}(\overline{a_{i}} + \overline{c_{i}ab}))\geq d-1.$ Thus $(\overline{a_{2}'}, \ldots, \overline{a_{d}'})\in \Um_{d-1}(\overline{R})$, where $ \overline{a_{i}'} = \overline{a_{i}} + \overline{c_{i}ab}$ for $2\leq i\leq d.$ Now as done in the previous paragraph, we can modify $v, w$ to $v_{3} =  ({a}, {a_{2}''}, \ldots, a_{d-1}'', {a_{d}''}), w_{3} = ({a}, {a_{2}''}, \ldots, a_{d-1}'', {a_{d}''})$ respectively such that $v_{3}, w_{3} \equiv e_{d}~\mbox{mod}~I$ and $(v,w)\rightarrow (v_{3}, w_{3})$ is an admissible sequence. $\hfill\square$

\begin{theo}\label{2}
\label{mainthm} Let $R$ be an affine algebra of dimension $d\geq 4$ over an algebraically closed field $k$ of char $\neq 2,3.$ Assume that either (i) $k = \overline{\mathbb{F}}_{p}$ or (ii) $R$ is normal. Then the group structure on ${\Um_{d}(R)}/{\E_{d}(R)}$ is nice.
\end{theo}

${\pf}$  In view of \cite[Lemma 3.5]{gr}, we may assume that $R$ is a reduced. Let $v, w \in \Um_{d}(R)$ and $I$ denotes the ideal defining the singular locus of $R.$ By (\ref{mainprop}), we may assume $(v, w)$ to be an admissible pair such that $v, w \equiv e_{d}~\mbox{mod}~I.$ Let $v = (a, a_{1}, \ldots, a_{d-1}), w = (b, a_{1}, \ldots, a_{d-1}).$  Using {\cite[Theorem 6.8]{ggr}} where the result is proved in the smooth case and (\ref{keylemma}) with $r(d) =d-1$ we get 
$ [v]*[w]= [(a, a_{1}, \ldots, a_{d-1})]\ast  [(b, a_{1}, \ldots, a_{d-1})] =  [(ab, a_{1}, \ldots, a_{d-1})].$ 
 $\hfill\square$

%\section{Niceness of relative unimodular space group}

\begin{theo}\label{relativemain}\label{21}
Let $R$ be an affine algebra of dimension $d\geq 4$ over $\overline{\mathbb{F}}_{p}$ with $p \neq 2,3$ and $I$ an ideal of $R.$ Then the group structure on ${\Um_{d}(R, I)}/{\E_{d}(R, I)}$ is nice. 
\end{theo}
${\pf}$ Since $R\oplus I$ is an affine $\overline{\mathbb{F}}_{p}$-algebra of dimension $d$ \cite[Proposition 3.1]{K}, it has nice group structure by (\ref{mainthm}). Now the result follows from (\ref{relnicecri}.)
$\hfill\square$

\begin{rem} Let $R$ be a normal affine algebra of dimension $d\geq 4$ over an algebraically closed field  $k$ of char $\neq 2,3.$ If $I$ is an ideal of $R$, then $R\oplus I$ is not normal and so we can not use (\ref{mainthm}) to conclude that $\Um_{d}(R, I)/\E_{d}(R, I)$ has nice group structure. Below we give a separate proof in this case.
\end{rem}

\begin{lem}
\label{relelementary} Let $R$ be a ring of dimension $d\geq 4$ and $J$ be an ideal of $R$ such that $\mbox{dim}(R/J)\leq d-2.$ Let $v\in \Um_{d}(R, I)$ for some ideal $I \subset R.$ Then there exists $\varepsilon \in \E_{d}(R, I)$ such that 
$v\varepsilon = (u_{1}, u_{2}, \ldots, u_{d})$ with $u_{1} \equiv 1~\mbox{mod}~(I \cap J).$
\end{lem}
${\pf}$ Let $v = (a_{1}, a_{2},\ldots, a_{d})\in \Um_{d}(R, I)$ and $\overline{R} = R/J.$ We have $\overline{v} = (\overline{a_{1}}, \overline{a_{2}},\ldots, \overline{a_{d}})\in \Um_{d}(\overline{R}, \overline{I}).$ By stable range condition \cite{bass}, upon adding multiples of $\overline{a_{d}}$ to $\overline{a_{1}}, \overline{a_{2}},\ldots, \overline{a_{d-1}}$,  we may assume that $(\overline{a_{1}}, \overline{a_{2}},\ldots, \overline{a_{d-1}})\in \Um_{d-1}(\overline{R}, \overline{I}).$ Let $(\overline{b_{1}}, \overline{b_{2}},\ldots, \overline{b_{d-1}})\in \Um_{d-1}(\overline{R}, \overline{I})$ be such that 
$\sum_{i=1}^{d-1} \overline{a_{i}} \overline{b_{i}} = \overline{1}.$ Let $$\overline{\varepsilon_{1}} = \begin{bmatrix}
                 \overline{1} &  \overline{0} & \ldots & (1-\overline{a_{d}})\overline{b_{1}}\\
                  \overline{0} & \overline{1}& \ldots & (1-\overline{a_{d}})\overline{b_{2}}\\
\vdots & & & \vdots \\
\overline{0} & \overline{0}& \ldots & (1-\overline{a_{d}})\overline{b_{d-1}}\\
\overline{0} & \overline{0}& \ldots & \overline{1}\\
                \end{bmatrix} \in \E_{d}(\overline{R}).$$ Then $(\overline{a_{1}}, \overline{a_{2}},\ldots, \overline{a_{d}})\overline{\varepsilon_1} = (\overline{a_{1}}, \overline{a_{2}},\ldots, \overline{a_{d-1}}, \overline{1}).$ There exists $\overline{\varepsilon_2}
\in\E_{d}(\overline{R}, \overline{I})$ such that $(\overline{a_{1}}, \overline{a_{2}},\ldots, \overline{a_{d-1}}, \overline{1})\overline{\varepsilon_2} = (\overline{1}, \overline{0},\ldots, \overline{0}, \overline{1}).$ Now $(\overline{1}, \overline{0},\ldots, \overline{0}, \overline{1})\overline{\varepsilon_1}^{-1} = (\overline{1}, \overline{0},\ldots, \overline{0}, \overline{1}-(1-\overline{a_{d}})\overline{b_1}).$ Since 
$1 -(1-\overline{a_{d}})\overline{b_1}) = \overline{0}~\mbox{mod}~{\overline{I}},$ there exists $\overline{\varepsilon_3}\in \E_{d}(\overline{R}, \overline{I})$ such that   $ (\overline{1}, \overline{0},\ldots, \overline{0}, 1-(1-\overline{a_{d}})\overline{b_1})\overline{\varepsilon_3} = (\overline{1}, \overline{0},\ldots, \overline{0}, \overline{0}).$  Let $\overline{\varepsilon} = \overline{\varepsilon_1}\overline{\varepsilon_2}\overline{\varepsilon_1}^{-1} \overline{\varepsilon_3} \in \E_{d}(\overline{R}, \overline{I}).$ Let  $\varepsilon \in \E_{d}(R, I)$ be a lift of $\overline{\varepsilon}.$ Then we have $v\varepsilon = (u_{1}, u_{2}, \ldots, u_{d})$ with $u_{1} \equiv 1~\mbox{mod}~(I \cap J).$
$~~~~~~~~~~~~~~~~~~~~~~~~~~~~~~~~~~~~~~~~~~~~~~~~~~~~~~~~~~~~~~~~~~~~~~~~~~~~~~~~~~~~~~~~~~~~~~~~~~~~~~~~~
           ~~~~~\qedwhite$

\begin{theo}\label{22}  
Let $R$ be a normal affine algebra of dimension $d\geq 4$ over an algebraically closed field $k$ of char $\neq 2,3$ and $I$ be an ideal of $R.$ Then the group structure on $\Um_{d}(R, I)/\E_{d}(R, I)$ is nice.
\end{theo} 
${\pf}$ Let $v = (a_{1}, a_{2}, \ldots, a_{d-1}, a), w = (a_{1}, a_{2}, \ldots, a_{d-1}, b)$ be such that $v, w \in \Um_{d}(R, I).$ Then $(a_{1}, a_{2}, \ldots, a_{d-1}, ab)\in \Um_{d}(R, I).$ Let $J$ be the ideal defining singular locus of $R.$ Since $R$ is a normal affine algebra, $\mbox{ht}(J)\geq 2.$  By (\ref{relelementary}), we may assume that $a_{1} \equiv 1~\mbox{mod}~(I \cap J).$  Let $a_{1} = 1 - \lambda$ for some $\lambda \in I\cap J.$ Then $(a_{1}, \lambda a_{2}, \ldots, \lambda a_{d-1}, \lambda ab)\in 
\Um_{d}(R, I).$ By Swan's Bertini (\ref{sb}), we can add a linear combination of $\lambda ab, \lambda a_{d-1}, \ldots, \lambda a_{2}$ to $a_{1}$ which transforms it to $a_{1}'$ such that the hypersurface $\ol R=R/\langle a_{1}'\rangle$ is smooth outside the singular set of $R$ and $\mbox{dim}(\ol R) = d-1.$ Since $a_{1}' \equiv 1~\mbox{mod}~ J,$ $\ol R$ is a smooth affine algebra. Let $a_{1}' = 1-\eta.$ Note that $\overline{\eta} = \overline{1}.$ Now choose $(b_{1}, b_{2}, \ldots, b_{d-1}, p)\in \Um_{d}(R, I)$ such that $a_{1}'b_{1} + \sum_{i = 2}^{d-1}a_{i}b_{i} + ap = 1.$ Then 
\begin{eqnarray*} 
 [w] \ast [v] &=& [(a_{1}', a_{2}, \ldots, a_{d-1}, b)] \ast [(a_{1}', a_{2}, \ldots, a_{d-1}, a)] \\
 &=& [(a_{1}', a_{2}, \ldots, a_{d-2}, a_{d-1}(b+p), a(b+p)-\eta)] 
\end{eqnarray*} 
Upon going modulo $\langle a_{1}'\rangle$, we have 
$$
[( \overline{a_{2}}, \ldots, \overline{a_{d-1}}(\overline{b}+\overline{p}), \overline{a}(\overline{b}+\overline{p})-\overline{\eta})]   = [(\overline{a_{2}}, \ldots, \overline{a_{d-1}}, \overline{b})] \ast [(
\overline{ a_{2}}, \ldots, \overline{a_{d-1}}, \overline{a})] = [(\overline{a_{2}}, \ldots, \overline{a_{d-1}}, \overline{ab})]
$$
by {\cite[Theorem 6.8]{ggr}}.
Since $\overline{\eta} = \overline{1},$ we have $$(a_{2}, \ldots, a_{d-2}, a_{d-1}(b+p), a(b+p)-\eta)\alpha \equiv ( a_{2}, \ldots, a_{d-1}, ab)~~\mbox{mod}~(Ia_{1}')$$ for some $\alpha \in \E_{d-1}(\eta) \subset \E_{d-1}(I).$ Thus 
$$ [(a_{1}', a_{2}, \ldots, a_{d-2}, a_{d-1}(b+p), a(b+p)-\eta)]  = [(a_{1}', a_{2}, \ldots, a_{d-1}, ab)]
 = [(a_{1}, a_{2}, \ldots, a_{d-1}, ab)]. $$
$~~~~~~~~~~~~~~~~~~~~~~~~~~~~~~~~~~~~~~~~~~~~~~~~~~~~~~~~~~~~~~~~~~~~~~~~~~~~~~~~~~~~~~~~~~~~~~~~~~~~~~~~~
           ~~~~~\qedwhite$

 \begin{cor}
           \label{vandercor}
            Let $R$ be an affine algebra of dimension $d\geq 4$ over an algebraically closed field $k$ of char $\neq 2, 3.$ Let $I$ be an ideal of $R.$ Assume that either (i) $k = \overline{\mathbb{F}}_{p}$ or (ii) $R$ is normal. If $\sigma \in {\SL}_{d}(R, I)\cap {\E}_{d+1}(R, I)$, then 
          $e_{1}\sigma$ can be completed to a relative elementary matrix $\varepsilon \in \E_{d}(R, I).$
           \end{cor}
    ${\pf}$ Proof is similar to (\ref{vandercor1}).
$\hfill \square$

\section{The first row map}

\iffalse

Let $R$ be a ring of dimension $d\geq 2$.
In \cite[Theorem 3.16(iv)]{vdk1}, van der Kallen proved that the 
first row map 
$\SL_{d+1}(R)\longrightarrow {\Um_{d+1}(R)}/{\E_{d+1}(R)}$
is a group homomorphism. Further, he gave an example showing that the first row map 
$\GL_{d+1}(R)\longrightarrow {\Um_{d+1}(R)}/{\E_{d+1}(R)}$
is not a group homomorphism \cite[Examples 4.16, 4.13]{vdk1}.

\par In \cite[Theorem 1.1]{grouphom},  the second author proved that if $R$ is local, then the first row map $\GL_{d+1}(R[X]) \rightarrow \Um_{d+1}(R[X])/{\E_{d+1}(R[X])}$ is a group homomorphism. Using this result we prove the following : 
\fi
%\subsection{Group homomorphism $row_1$}

\begin{theo}
\label{gphomcors} Let $R$ be a ring of dimension $d\geq 2$ such that $\E_{d+1}(R)$ acts transitively on $\Um_{d+1}(R).$ Then the first row map 
$row_1:\GL_{d+1}(R[X]) \rightarrow \Um_{d+1}(R[X])/{\E_{d+1}(R[X])}$ is a group homomorphism.
\end{theo}
${\pf}$ Let $\sigma, \tau \in \GL_{d+1}(R[X])$ and $[v] = [e_{1}\sigma]\ast [e_{1}\tau] \ast [e_{1}\sigma\tau]^{-1}.$ For $\mathfrak{m}\in Max(R),$ by \cite[Theorem 1.1]{grouphom} where the result is proved when $R$ is local, $[v]_{\mathfrak{m}} = [e_{1}].$ Therefore by (\ref{raolg}),
$[v] = [v({0})].$ Since $\E_{d+1}(R)$ acts transitively on $\Um_{d+1}(R),$ we have $[v] = [e_1].$ Thus $[e_{1}\sigma]\ast [e_{1}\tau] = [e_{1}\sigma\tau].$ $~~~~~~~~~~~~~~~~~~~~~~~~~~~~~~~~~~~~~~~~~~~~~~~~~~~~~~~~~~~~~~~~~~~~~~~~~~~~~~~~~~~~~~~~~~~~~~~~~~~~~~~~~
           ~~~~~\qedwhite$

\begin{cor}\label{elem3}
Let $R$ be a ring of dimension $d\geq 2$. Assume that either (i) $R$ is an affine $C$-algebra with $C=\mathbb{ Z}$ or $C$ is a subfield of $\overline{\mathbb{F}}_{p}$  or (ii) Jacobson radical has $\mbox{ht}(J(R))\geq 1$. Then $\E_{d+1}(R)$ acts transitively on $\Um_{d+1}(R).$ As a consequence, the first row map 
$row_1:\GL_{d+1}(R[X]) \rightarrow \Um_{d+1}(R[X])/{\E_{d+1}(R[X])}$ is a group homomorphism.
\end{cor}
${\pf}$  Using \cite[Corollary 18.1, Theorem 18.2]{7} for $C = \mathbb{Z}$, \cite[Corollary 17.3]{7} when $C$ is a subfield of $\overline{\mathbb{F}}_{p}$ and \cite[Lemma 4.2]{kesharisharma} when $\mbox{ht}(J(R))\geq 1$, we get  $\E_{d+1}(R)$ acts transitively on $\Um_{d+1}(R).$
$\hfill\square$

\section{Torsion in $\Um_{d+1}(R[X])/{\E_{d+1}(R[X])}$}

\begin{theo}
\label{torcors} Let $R$ be a ring of dimension $d\geq 2$ with $1/m,1/d!\in R$ such that $\E_{d+1}(R)$ acts transitively on $\Um_{d+1}(R).$ Then the group ${\Um_{d+1}(R[X])}/{\E_{d+1}(R[X])}$ has no $m$-torsion.
\end{theo}

${\pf}$ Let $[v]\in {\Um_{d+1}(R[X])}/{\E_{d+1}(R[X])}$ such that $[v]^{m} = [e_{1}].$ For $\mathfrak{m}\in Max(R),$ by \cite[Corollary 3.8]{absence} where the result is proved when $R$ is local, $[v]_{\mathfrak{m}} = [e_{1}].$ Therefore by (\ref{raolg}),
$[v] = [v({0})].$ Since $\E_{d+1}(R)$ acts transitively on $\Um_{d+1}(R),$ we have $[v] = [e_1].$  $~~~~~~~~~~~~~~~~~~~~~~~~~~~~~~~~~~~~~~~~~~~~~~~~~~~~~~~~~~~~~~~~~~~~~~~~~~~~~~~~~~~~~~~~~~~~~~~~~~~~~~~~~
           ~~~~~\qedwhite$

\begin{cor}\label{3}
\label{gargetype}  Let $R$ be a ring of dimension $d\geq 2$ with $1/m,1/d!\in R$. Assume that either (i) $R$ is an affine $C$-algebra with $C=\mathbb{ Z}$ or $C$ is a subfield of $\overline{\mathbb{F}}_{p}$  or (ii) Jacobson radical has $\mbox{ht}(J(R))\geq 1$. Then 
the group ${\Um_{d+1}(R[X])}/{\E_{d+1}(R[X])}$ has no $m$-torsion.

\end{cor}
${\pf}$  Use (\ref{elem3}, \ref{torcors}).
$\hfill \square$

\section{Nice group structure on ${\Um_{d+1}(R[X], I[X])}/{\E_{d+1}(R[X], I[X])}$}

%The absolute case of next result is proved in {\cite[Corollary 4.5]{kesharisharma}}.

\begin{theo}\label{31}
Let $R$ be a ring of dimension $d\geq 4$ with $1/2\in R$ and $I$ be an ideal of $R$. Assume that either (i) $R$ is an affine $C$-algebra with $C=\mathbb{ Z}$ or $C$ is a subfield  of $\overline{\mathbb{F}}_{p}$  or (ii) Jacobson radical has $\mbox{ht}(J(R))\geq 1$.
Then the group structure on $\Um_{d+1}(R[X],I[X])/{\E}_{d+1}(R[X],I[X])$ is nice.
\end{theo}
${\pf}$ When $R$ is an affine $C$-algebra of dimension $d$, then $R\oplus I$ is also an affine $C$-algebra of dimension $d$. Further $R\oplus I$ is an integral extension of $R$ \cite[Proposition 3.1]{K}. Therefore the group structures on 
$\Um_{d+1}(R[X])/{\E}_{d+1}(R[X])$ and
$\Um_{d+1}(R[X] \oplus I[X])/{\E}_{d+1}(R[X] \oplus I[X])$ are nice {\cite[Corollary 4.5]{kesharisharma}}. Now use 
(\ref{relnicecri}) to conclude the proof.
$\hfill \square$
\medskip

{\bf Acknowledgement:} We thank the referee for detailed comments some of which are included in the introduction and also for suggesting (\ref{keylemma}). The second author thanks DST inspire DST/INSPIRE/04/2021/002849 and SERB SRG/2022/000056 for their financial support.

\Addresses

\end{document}